\def\pf(#1,#2){\langle\langle #1,#2\rangle\rangle}
\def\pf(#1){\langle\langle #1\rangle\rangle}
\def\symb(#1){\langle #1\rangle}

\def\Eq{\mathop{\sl Eq}}

\def\tr{\mathop{\rm tr}}
\def\stufe(#1){s(#1)}

\def\G{SL(2,K)}
\def\PG{PSL(2,K)}
\def\MG{St(2,K)}
\def\Z{{\bf Z}}
\def\R{{\bf R}}

\def\TccI{{TccI}}
\def\TccII{{TccII}}
\def\dotK{{K^*}}
\def\dotKK{K^{*2}}
\def\fund#1{\pi_1(SL(2,#1))}
\def\ov{\overline}
\def\Q{{\bf Q}}

\def\qed{\hfill{$\diamond$}}

\def\c1c{\mathop\tau}
\def\w{\mathop{\sl w}}
\def\c{\mathop{\sl c}}
\def\u{u}

\def\com(#1){#1^{[,]}}

\def\fund{H_2SL(2,K)}
\def\pfund{H_2PSL(2,K)} 
\def\MG{St(2,K)}
\def\ov{\overline}

\def\wt{\widetilde}
\def\Q{{\bf Q}}

\newcount\m
\def\t{\the\m\global\advance\m by 1}
\m=1
\newcount\n
\def\f{\the\n\global\advance\n by 1}
\n=1

\centerline{\bf Tautological characteristic classes III:
the Witt class for PSL(2)}
\bigskip
\rm
Jan Dymara

Instytut Matematyczny Uniwersytetu Wroc\l awskiego,

pl.~Grunwaldzki 2, 50-384 Wroc\l aw

{\tt jan.dymara@uwr.edu.pl}
\smallskip

Tadeusz Januszkiewicz

Instytut Matematyczny PAN,

ul.~Kopernika 18, 51-617 Wroc\l aw

{\tt tjan@impan.pl}


\footnote{}{2020 \it Mathematics Subject Classification.
\rm Primary 20G10; Secondary 11E81. }
\footnote{}{\it Key words and phrases. \rm Witt class, equicommutative class,
central extension, Milnor--Wood inequality.}

\bigskip\rm
\bigskip\rm
\medskip
\hskip1cm\vbox{\hsize13.5cm
\noindent ABSTRACT. We explain the relation between the Witt class and
the universal equicommutative class for $PSL(2,K)$. We discuss
an analogue of the Milnor--Wood inequality.
}
\bigskip
\bigskip
\centerline{\bf Introduction}
\medskip
The Witt class of $\PG$ is a second cohomology class of this group,
with coefficients in the Witt ring of symmetric bilinear forms over $K$,
tautologically arising from the action of $\PG$ on the projective
line ${\bf P}^1(K)$ (cf [Ne], [Ba], [Kr-T], [\TccII]).
In [\TccII] we investigated its pull-back to $\G$.
It turns out that some properties studied there
do not always hold for $\PG$. The differences are significant
and are discussed in detail in the present paper.
Thus:
\item{$\bullet$} we calculate the coefficient group of the universal
equicommutative class of $\PG$ (see Corollary 3.6, Theorems 4.3, 4.5);
\item{$\bullet$} we explain the relation between the Witt class
and the universal equicommutative class for $\PG$, which turns out
to be subtly dependent on the arithmetic of $K$ (see Theorem 5.4);
\item{$\bullet$} we discuss the range of values of the Witt class
for $\PG$-bundles over a surface of genus $g$, establishing an analogue
of the Milnor--Wood inequality (see Theorem 6.1).
\smallskip
It is likely that the Witt class of $\PG$ will be a valuable tool
in the study of varieties of projective representations
over general fields $K$,
just as the Euler class is for $K=\R$.
\smallskip
Throughout the paper $K$ denotes an infinite field of characteristic
different from $2$. (For finite fields the Witt class is virtually trivial;
in characteristic $2$ we have $\G=\PG$, so that the results obtained
in [\TccII] for $\G$ do not require any adjustments.)

\bigskip

\medskip
\bf 1. Universal extensions and classes.\rm
\medskip
\def\a{1}
\m=1
\n=1

A central group extension $\pi\colon E\to G$ is called \it equicommutative \rm
if
$$\forall g,h\in E,\quad\pi(g)\pi(h)=\pi(h)\pi(g)\,\Rightarrow\,gh=hg.\leqno(\a.\f)$$
Equicommutative extensions are important because they correspond to
twist--invariant characteristic classes of flat $G$-bundles (cf [\TccII, Section 5]).
A perfect group $G$ has a universal central extension $\pi\colon\overline{G}\to G$,
with kernel $H_2G$ ($=H_2(G,\Z)$) called the Schur multiplier of $G$.
This extension is, in general, not equicommutative; it has, however, a canonical
largest equicommutative quotient. To obtain this quotient one divides $\overline{G}$
by a certain subgroup $\com(G)$ of $\ker{\pi}=H_2G$.
One description of $\com(G)$ is that it is generated by the following set of commutators:
$$\{[a,b]\mid a,b\in \overline{G}, [\pi(a),\pi(b)]=1\}.\leqno(\a.\f)$$
Once we fix a set--theoretic section (lift) $G\ni g\mapsto\ov{g}\in\ov{G}$,
 we can give a more economical formula for the set (\a.2), namely
$$\{[\ov{g},\ov{h}]\mid g,h\in G, [g,h]=1\}.\leqno(\a.\f)$$
Indeed, if $\pi(a)=\pi(a')$ and $\pi(b)=\pi(b')$, then, using the assumption
that the extension is central, one shows $[a,b]=[a',b']$.
Another description of $\com(G)$: it is the subgroup of $H_2G$ generated by
toral classes, ie by images of the fundamental class of the
$2$-dimensional torus
$T^2$ in $H_2G=H_2(BG,\Z)$ under all maps $T^2\to BG$. (Cf [\TccII, Definition 6.1, Remark 6.3].)

To summarize, for a perfect $G$ the extension $\overline{G}/\com(G)\to G$
is the universal equicommutative central extension; its kernel is equal to $H_2G/\com(G)$
and is denoted $\Eq(G)$.

The above discussion of central extensions has a cohomological counterpart (cf [Brown, Chapter IV,\S 3]).
Let $\ov{G}\to G$ be a central extension with kernel $U$; choose a lift $g\mapsto\ov{g}$
and define the cocycle $c$ by the formula
$$\ov{g}\cdot\ov{h}=\ov{gh}\cdot c(g,h).\leqno(\a.\f)$$
The cohomology class $[c]\in H^2(G,U)$ does not depend on the choice of the lift.
This construction provides a bijective correspondence between
isomorphism classes of central extensions of $G$ with kernel $U$ and cohomology classes in $H^2(G,U)$.
Equicommutative extensions correspond to equicommutative classes, ie classes represented by
equicommutative cocycles: a cocycle $c$ is equicommutative if
$$\forall g,h\in G\quad gh=hg\Rightarrow c(g,h)=c(h,g).\leqno(\a.\f)$$
For a perfect group $G$, the universal central extension corresponds to a class $u_G\in H^2(G,H_2G)$
that enjoys the following universality property:
for any cohomology class $u\in H^2(G,U)$ there exists a unique coefficient group
homomorphism $H_2G\to U$ that maps $u_G$ to $u$.
We call $u_G$ the universal class or the Moore class of $G$.
Similarly, the universal equicommutative central extension corresponds to
a class $w_G\in H^2(G,\Eq(G))$ that has a universality property analogous to that of $u_G$,
but for equicommutative classes $u$. We call $w_G$ the universal equicommutative class of $G$.
More details, proofs and references may be found in [\TccII, Sections 5 and 6].

\bigskip

\medskip
\bf 2. $H_2SL(2,K)$ and $H_2PSL(2,K)$.\rm
\medskip
\def\a{2}
\m=1
\n=1

Let
$K$ be an infinite field of characteristic different from $2$.
In this section we recall the classical description of the universal central
extension and of the Schur multiplier of $\G$ and of $\PG$.
The standard references are
[Moore, Sections 8,9], [Mats], [St, \S 7]. For the reader's convenience, we repeat the
formulae from [TccII] (with adjusted numbering and references):

\bigskip
\hskip.5cm\vbox{\hsize14.5cm
The universal central extension of $SL(2,K)$ is called the Steinberg group
and is denoted $St(2,K)$. It is generated by two families of symbols:
$x_{12}(t), t\in K$; $x_{21}(t), t\in K$. For $t\in\dotK$ one defines
an auxiliary element $w_{ij}(t)=x_{ij}(t)x_{ji}(-t^{-1})x_{ij}(t)$; then
the relations defining $St(2,K)$ are:
$$\eqalign{
x_{ij}(t)x_{ij}(s)&=x_{ij}(t+s)\qquad(t,s\in K),\cr
w_{ij}(t)x_{ij}(r)w_{ij}(t)^{-1}&=x_{ji}(-t^{-2}r)\qquad(t\in\dotK,r\in K),}
\leqno(\a.\f)$$
where $\{i,j\}=\{1,2\}$. Other noteworthy elements of $St(2,K)$ are
$h_{ij}(t)=w_{ij}(t)w_{ij}(-1)$.

The projection $\pi\colon St(2,K)\to SL(2,K)$ is defined by
$$\pi\colon\quad x_{12}(t)\mapsto\pmatrix{1&t\cr0&1},
\quad x_{21}(t)\mapsto\pmatrix{1&0\cr t&1}.\leqno(\a.\f)$$
Then one easily checks that
$$\pi\colon\quad w_{12}(t)\mapsto\pmatrix{0&t\cr-t^{-1}&0},
\quad h_{12}(t)\mapsto\pmatrix{t&0\cr 0&t^{-1}}.\leqno(\a.\f)$$
}

\hskip.5cm\vbox{\hsize14.5cm
The kernel of $\pi$, i.e.~the Schur multiplier of $SL(2,K)$
(denoted $H_2(SL(2,K))$, $\pi_1(SL(2,K))$ or $KSp_2(K)$ in various sources)
is generated by elements
$$\{s,t\}:=h_{12}(s)h_{12}(t)h_{12}(st)^{-1}\quad (s,t\in\dotK).\leqno(\a.\f) $$
With this generating set it is described by an explicit family of relations:
$$\eqalign{
&\{st,r\}\{s,t\}=\{s,tr\}\{t,r\},\quad \{1,s\}=\{s,1\}=1;\cr
&\{s,t\}=\{t^{-1},s\};\cr
&\{s,t\}=\{s,-st\};\cr
&\{s,t\}=\{s,(1-s)t\} \quad \hbox{\rm if $s\ne1$}.}\leqno(\a.\f)
$$
\noindent (For this presentation see [Moore, Theorem 9.2] or
[St, \S 7, Theorem 12]. Even though the group is abelian, the
convention is multiplicative.)\hskip.2cm([\TccII, Section 7])
}
\bigskip

Let $p\colon\G\to\PG$ be the natural projection. The composition
$p\circ\pi\colon \MG\to\PG$ is the universal central extension [St, \S 7 (vi)].
Since $\PG$ has trivial centre, the group $\pfund=\ker(p\circ\pi)$
is equal to the centre of $\MG$. We have a short exact sequence
of abelian groups:
$$0\to\fund\to\pfund\to\Z/2\to0.\leqno(\a.\f)$$
We can describe $\pfund$ as a central extension of $\Z/2$ by $\fund$
in term of a cocycle. To do this we lift both elements
of $\Z/2$ to $\pfund$: $+I$ to $1$, and $-I$ to $h_{12}(-1)$.
(We shall often abbreviate $h_{12}(t)$ to $h(t)$; in particular $h_{12}(-1)$ will be $h(-1)$.)
Then the extension is described by the cocycle on $\Z/2$ with
$\fund$--coefficients whose only non-trivial value is taken on
$(-I,-I)$ and is equal to
$$h(-1)h(-1)h((-1)\cdot(-1))^{-1}=\{-1,-1\}.\leqno(\a.\f)$$
Thus $\pfund$ is an abelian group generated by the symbols
$\{s,t\}$ and the additional
generator $h(-1)$. The relations are
the symbol relations (\a.5) (as in $\fund$) and
the extra relation $h(-1)^2=\{-1,-1\}$ (note that $h(1)=1$).

It is tempting to ask about the nature of the sequence (\a.6), eg if it splits.
One easily checks that it does split if and only if the symbol
$\{-1,-1\}$ is a square in the group $\fund$.
We expect that the latter condition depends on
arithmetic properties of the field $K$.
For example, if $-1$ is a square in $\dotK$, say $-1=a^2$, then, by
[Moore, Lemma 3.2], $\{-1,-1\}=\{a^2,a^2\}=\{a^2,a\}^2$. 
\medskip

We have already mentioned that the Steinberg group $\MG$ is the universal central extension
both of $\G$ and of $\PG$. We shall encounter other central extensions of $\G$ that it will
be necessary to consider as extensions of $\PG$. The following remark will then be useful.
 
\medskip

\bf Remark \a.\t.\rm\par
Let $\pi\colon G\to\G$ be a central extension; then $p\circ\pi\colon G\to \PG$ is also
a central extension. This follows from [Wei, Exercise 6.9.2]. (The correct version of
the exercise assumes additionally that $X$---corresponding to our $\G$---is perfect.)

\bigskip

\medskip
\bf 3. $\Eq(SL(2,K))$ and $\Eq(PSL(2,K))$.\rm
\medskip
\def\a{3}
\m=1
\n=1
In this section we begin to compute $\Eq(\PG)$---the kernel
of the universal equicommutative extension of $\PG$, ie the coefficient group of
the universal equicommutative cohomology class (these notions were explained in Section 1).

The group $\Eq(\G)$ has been calculated in [\TccII].
We recall its two descriptions, cf [\TccII, Theorem 10.1, Proposition 8.1]. First, it is a quotient
of $H_2\G$---the group generated by symbols $\{s,t\}$ with relations
(2.5)---obtained by adding extra relations of the form $\{s,t\}=\{t,s\}$.
Second, it is isomorphic to $I^2(K)$, the square of the fundamental ideal
$I(K)$ of the Witt ring $W(K)$ of the field $K$. The isomorphism
is induced by the map $\{s,t\}\mapsto\pf(s,t)$, where $\pf(s,t)$ denotes
the Pfister form---the 4-dimensional diagonal quadratic form
with coefficients $\langle 1,-s,-t,st\rangle$.

We would like to obtain
some analogue of the second description for $\Eq(\PG)$.
To this end, we pass
to the quotient by $\G^{[,]}$ with the first two non-trivial terms
of the sequence (2.6);
we obtain
$$0\to\Eq(\G)\to\pfund/\G^{[,]}\to\Z/2\to0.\leqno(\a.\f)$$
Since $\G^{[,]}\subseteq\PG^{[,]}$, the group $\Eq(\PG)=\PG/\PG^{[,]}$
is a quotient of the middle term $\pfund/\G^{[,]}$.
To proceed, it is crucial to find out whether $\G^{[,]}=\PG^{[,]}$.
This turns out to depend on the arithmetic of the field $K$.

In order to compare the groups $\G^{[,]}$ and $\PG^{[,]}$
we shall now describe them in terms of lifts.
We denote elements of $\PG$ by $\pm A$, where $A\in \G$.
The neutral element of $\PG$ is $\pm I$. Elements
$\pm A$, $\pm B$ commute if and only if $AB=BA$ or $AB=-BA$;
notice that the sign in this equation does not change when we replace
$A$ by $-A$, or $B$ by $-B$.
Let us fix a set-theoretic lift $\G\ni A\mapsto\ov{A}\in \MG$.
Then, since the lift $h(-1)$ of $-I$ is central in $\MG$,
the lifts $\ov{A}$ and $\ov{-A}$ differ by a central element.
It follows that the four commutators
$$[\ov{A},\ov{B}],\quad[\ov{-A},\ov{B}],\quad[\ov{A},\ov{-B}],\quad
[\ov{-A},\ov{-B}]$$
are equal. Let us pick an (arbitrary) function $\sigma\colon\G\to\Z/2=\{1,-1\}$
that satisfies $\sigma(I)=1$ and $\sigma(-A)=-\sigma(A)$ for all $A\in\G$.
Then $\pm A\mapsto \sigma(A)A$ is a well-defined lift $\PG\to \G$,
and $\pm A\to \ov{\pm A}:=\ov{\sigma(A) A}$ is a lift $\PG\to \MG$.
We can now see that 
$$\eqalign{
\com(\G)&=\langle\{[\ov{A},\ov{B}]\mid A,B\in \G, AB=BA\}\rangle\cr
&=\langle\{[\ov{\pm A},\ov{\pm B}]\mid \pm A,\pm B\in \PG, AB=BA\}\rangle.}
\leqno(\a.\f)$$
To enlarge this group to $\PG^{[,]}$ we need to add the
commutators $[\ov{\pm A},\ov{\pm B}]$ for pairs $A$, $B$
that satisfy $AB=-BA$:
$$\eqalign{
\com(\PG)&=\langle\com(\G)\cup\{[\ov{A},\ov{B}]\mid A,B\in \G, AB=-BA\}
\rangle\cr
&=\langle\com(\G)\cup\{[\ov{\pm A},\ov{\pm B}]\mid \pm A,\pm B\in \PG, AB=-BA\}
\rangle.}
\leqno(\a.\f)$$

\smallskip
Do anti-commuting pairs exist?

\smallskip
\bf Lemma \a.\t.
\sl\par\nobreak
The following are equivalent:
\item{1)} There exist $a,b\in SL(2,K)$ such that $ab=-ba$;
\item{2)} $-1$ is a sum of two squares in $K$.
\medskip
\rm Proof.\par
$2)\Rightarrow 1)$: Let $-1=\gamma^2+\delta^2$. 
Let
$$a=\pmatrix{\gamma&\delta\cr \delta&-\gamma},
\qquad b=\pmatrix{0&-1\cr1&0}.\leqno(\a.\f)$$ 
Then $ab=-ba$.

$1)\Rightarrow 2)$: 
Suppose $ab=-ba$, $a,b\in SL(2,K)$. We rewrite this equality as $aba^{-1}=-b$,
and we deduce $\tr{b}=-\tr{b}=0$. Hence the characteristic polynomial
of $b$ is $x^2+1$; the Cayley--Hamilton Theorem implies $b^2=-I$.
It follows that $b$ is $\G$-conjugate to
$$b'=\pmatrix{0&t\cr -t^{-1}&0}\leqno(\a.\f)$$
for some $t\in\dotK$; then the same conjugation transforms
$a$ to a matrix $a'$ that anti-commutes with $b'$.
Straightforward calculation shows that
$$a'=\pmatrix{\gamma&t^2\delta\cr \delta&-\gamma}\leqno(\a.\f)$$
for some $\gamma,\delta\in K$. The condition $\det{a'}=1$ is equivalent to 
$\gamma^2+(t\delta)^2=-1$.\qed
\smallskip
\bf Corollary \a.\t.\sl\par
Any pair of anti-commuting matrices in $\G$ is conjugate to a pair
given by $(\a.5)$, $(\a.6)$.
\rm\smallskip
The corollary follows from the proof of the implication $1)\Rightarrow 2)$
in the proof of Lemma \a.1.
\medskip
In view of Lemma \a.1, the following classical definition is obviously relevant.
\smallskip
\bf Definition \a.\t. \rm (cf.~[Lam, Definition XI.2.1])\par\nobreak
The \it Stufe \rm of a field $K$, denoted $\stufe(K)$, is the minimal length
of a representation of $-1$ as a sum of squares of elements of $K$.
If such representation does not exist, the Stufe is infinite.
\smallskip
It follows from Lemma \a.1 and from Definition \a.3 that for a field
$K$ of Stufe $>2$ there are no anti-commuting pairs of matrices in
$\G$. Consequently, in that case $\PG^{[,]}=\G^{[,]}$. This proves the
following proposition. 
\smallskip
\bf Proposition \a.\t.\sl\par
Let $K$ be an infinite field of 
Stufe $>2$.
Then $$\Eq(\PG)\simeq \pfund/\G^{[,]}.\leqno(\a.\f)$$
\smallskip
\rm
Recall that $\Eq(\G)$ is isomorphic to $I^2(K)$ via $\{s,t\}\mapsto\pf(s,t)$.
If $\stufe(K)>1$, i.e.~if $-1\not\in\dotK^2$, we extend this map to an injection
$\pfund/\G^{[,]}\to W(K)$ by sending
$h(-1)$ to $\langle 1,1\rangle$. This is indeed well--defined, since
$h(-1)^2=\{-1,-1\}\mapsto \pf(-1,-1)=2\langle 1,1\rangle$, and injective,
since $d_\pm(\langle 1,1\rangle)=-1$, hence
$\langle 1,1\rangle\not\in I^2(K)$ (here we use $\stufe(K)>1$).
We define $I^2_+(K)$ as the subgroup of $W(K)$ generated by
$I^2(K)$ and $\langle 1,1\rangle$. The signed determinant
restricts to a surjective map $d_{\pm}\colon I^2_+(K)\to\Z/2=\{1,-1\}$
with kernel $I^2(K)$.

The above discussion and Proposition \a.4 prove
the following proposition and corollary:
\smallskip
\bf Proposition \a.\t.\sl\par
Let $K$ be an infinite field of 
Stufe $>1$.
We have an isomorphism
of short exact sequences:
\def\mapdown#1{\Big\downarrow\rlap{$\vcenter{\hbox{$\scriptstyle{#1}$}}$}}
$$

\matrix{
0&\longrightarrow&
\Eq(\G)&\longrightarrow&\pfund/\G^{[,]}&\longrightarrow&\Z/2&\longrightarrow
&0\cr
{}&&\mapdown{}&&\mapdown{}&&\mapdown{\rm Id}&&
\cr
0&\longrightarrow&
I^2(K)&\longrightarrow&I^2_+(K)&
\longrightarrow&\Z/2&\longrightarrow&0\cr
}\leqno(\a.\f)
$$
\smallskip
\bf Corollary \a.\t.\sl\par
Let $K$ be an infinite field of 
Stufe $>2$.
Then $\Eq(\PG)\simeq I^2_+(K)$.\rm

\bigskip
Determining $\Eq(\PG)$ for $K$ of Stufe 1 or 2 requires extra arguments,
but even in these cases the first step is to
describe $\pfund/\G^{[,]}$. If $\stufe(K)=2$ this group is $I^2_+(K)$, by
Proposition \a.5. What happens for  $\stufe(K)=1$?
In that case $-1$ is a square and $\langle 1,1\rangle$
is $0$ in $W(K)$. Formally, we then have $I^2_+(K)=I^2(K)$.
The group $\pfund/\G^{[,]}$, however, is still an extension of $\Z/2$ by
$\Eq(\G)\simeq I^2(K)$, given by the 
cocycle with (the only potentially non-trivial) value $\pf(-1,-1)$
(the image of $\{-1,-1\}$ under
$\fund\to\Eq(\G)\simeq I^2(K)$).
But $-1\in\dotK^2$ implies $\pf(-1,-1)=0$ in $I^2(K)$, so that
the extension splits and $\pfund/\G^{[,]}\simeq I^2(K)\oplus\Z/2$.
The $\Z/2$ summand is generated by the image of $h(-1)$---a 2-torsion element.
Thus we have proved:
\smallskip
\bf Proposition \a.\t.\sl\par
Let $K$ be an infinite field of 
Stufe $1$.
Then $$\pfund/\G^{[,]}\simeq I^2(K)\oplus\Z/2.\leqno(\a.\f)$$
\rm\par
\medskip
\bf Remark \a.\t.\rm\par\nobreak
In is natural to inquire whether the lower row in (\a.8)
is a split extension. For $\stufe(K)=2$ it is indeed split.
In this case we have $-1=x^2+y^2$ for some $x,y\in K$. Then
a standard calculation (using the Witt relation) shows that $W(K)$ is a 4-torsion group:
$$
\eqalign{4\langle a\rangle&=
\langle a,a,a,a\rangle=
\langle a\rangle +\langle a \rangle - \langle -a\rangle - \langle -a\rangle\cr
&=\langle a\rangle +\langle a \rangle - \langle ax^2+ay^2\rangle - \langle ax^2ay^2(ax^2+ay^2)\rangle
=\langle a\rangle +\langle a \rangle - \langle ax^2\rangle - \langle ay^2\rangle\cr
&=\langle a\rangle +\langle a \rangle - \langle a\rangle - \langle a\rangle=0.}\leqno(\a.\f)
$$
It follows that
$2\langle1,1\rangle=4\langle1\rangle=0$, which implies the required splitting.
\bigskip

\medskip
\bf 4. $\Eq(PSL(2,K))$ for Stufe 2 and 1.\rm
\medskip
\def\a{4}
\m=1
\n=1

We assume throughout this section that $K$ is an infinite field
(of characteristic $\ne2$) and $\stufe(K)\le2$. Our goal is to describe $\Eq(\PG)=\pfund/\com(\PG)$ for such $K$. 
Formula (3.3) suggests a direct way to compute this group. We work in the universal extension
$\MG\to\G$. We compute $\com(\G)$; then we find the set of all commutators of lifts to $\MG$
of anti--commuting pairs in $\G$; we span a subgroup on this set and on $\com(\G)$---this will be $\com(\PG)$---and
then we are ready to describe $\pfund/\com(\PG)$.

We will follow an equivalent, simpler procedure. We work in the central extension
$$E=\MG/\com(\G)\to\G.$$ In this $E$ we have a subgroup $\pfund/\com(\G)$, identified by
Propositions 3.5, 3.7. The quotient map $\MG\to E$ maps:
commutators of lifts to $\MG$ of anti--commuting pairs in $\G$---to
commutators of lifts to $E$ of anti--commuting pairs in $\G$.
The latter commutators span the subgroup $\com(\PG)/\com(\G)$ of $E$;
the quotient of $\pfund/\com(\G)$ by this subgroup is $\pfund/\com(\PG)=\Eq(\PG)$.

We execute this procedure in three steps:
\item{1)} Prove two general Propositions \a.1, \a.2; these will facilitate
computations with lifts of anti--commuting pairs.
\item{2)} Recall from [\TccII] an explicit cocycle $\Phi_*b$ that corresponds to
the extension $E=\MG/\com(\G)\to\G$.
\item{3)} Use Corollary 3.2 and the cocycle $\Phi_*b$
to explicitly compute commutators of lifts to $E$ of all
anti-commuting pairs.
\bigskip
Step 1:
\smallskip
\bf Proposition \a.\t\sl\par\nobreak
Let a central extension  $E\to \G$ be given, with a lift
$\G\ni g\mapsto\ov{g}\in E$
and the corresponding cocycle $c$. Suppose $a,b\in \G$ satisfy $ab=-ba$.
Then
$$[\ov{a},\ov{b}]=\ov{-I}\,\cdot c(-I,ba)^{-1}c(a,b)c(b,a)^{-1}.\leqno(\a.\f)$$
\rm\par
Proof. We have 
$$\ov{a}\cdot\ov{b}=\ov{ab}\,c(a,b)=\ov{-ba}\,c(a,b)=
\ov{-I\cdot ba}\,c(a,b)
=\ov{-I}\cdot\ov{ba}\,c(-I,ba)^{-1}c(a,b).\leqno(\a.\f)$$
On the other hand
$\ov{b}\cdot\ov{a}=\ov{ba}\,c(b,a)$, 
hence
$$[\ov{a},\ov{b}]=\ov{a}\cdot\ov{b}\cdot(\ov{b}\cdot\ov{a})^{-1}=
\ov{-I}\cdot\ov{ba}\,c(-I,ba)^{-1}c(a,b)\cdot (\ov{ba}\,c(b,a))^{-1}=
\ov{-I}\cdot c(-I,ba)^{-1}c(a,b)c(b,a)^{-1}.\leqno(\a.\f)
$$
\qed
\smallskip
\bf Proposition \a.\t\sl\par\nobreak
Let a central extension $E\to \G$ be given,
with a lift $\G\ni g\mapsto\ov{g}\in E$
and the corresponding cocycle $c$. Suppose $a,b,d\in \G$, and $ab=ba$ or $ab=-ba$.
Then
$$[\ov{dad^{-1}},\ov{dbd^{-1}}]=[\ov{a},\ov{b}].\leqno(\a.\f)$$
\rm\par
Proof. The extension is central, hence the commutator
of lifts of two given elements does not depend on the choice of the lifts.
Notice that $\ov{d}\cdot\ov{a}\cdot\ov{d}^{-1}$ is a lift of $dad^{-1}$;
therefore
$$[\ov{dad^{-1}},\ov{dbd^{-1}}]=[\ov{d}\cdot\ov{a}\cdot\ov{d}^{-1},
\ov{d}\cdot\ov{b}\cdot\ov{d}^{-1}]=
\ov{d}[\ov{a},\ov{b}]\ov{d}^{-1}=[\ov{a},\ov{b}].\leqno(\a.\f)$$
The last equality holds because $[a,b]\in\{I,-I\}$ is central in $\G$, hence
its lift $[\ov{a},\ov{b}]$ is central in the extension (by Remark 2.1).
\qed
\medskip

Step 2:

We need an explicit form of a cocycle $b$
representing the class $\u_{\G}\in H^2(\G,\fund)$ corresponding
to the universal central extension $St(2,K)\to \G$---or
 rather the image of $b$ under the coefficient quotient map
$\Phi\colon\fund\to \fund/\G^{[,]}\simeq I^2(K)$.

Every element of $\G$ is uniquely represented in one of the forms:
$$g_1(u,t)=x(u)h(t)=\pmatrix{1&u\cr0&1}\pmatrix{t&0\cr0&t^{-1}};\leqno(\a.\f)$$
$$g_2(u,t,v)=x(u)w(t)x(v)=\pmatrix{1&u\cr0&1}\pmatrix{0&t\cr -t^{-1}&0}
\pmatrix{1&v\cr0&1}.\leqno(\a.\f)$$
This leads to the following 
lift $SL(2,K)\to St(2,K)$:
$$\ov{g_1(u,t)}=\ov{x(u)h(t)}:=x_{12}(u)h_{12}(t);\leqno(\a.\f)$$
$$\ov{g_2(u,t,v)}=\ov{x(u)w(t)x(v)}:=x_{12}(u)w_{12}(t)x_{12}(v)
\leqno(\a.\f)$$
The corresponding cocycle $b$ was calculated by Moore (see [Moore, 9.1-4]),
with later correction
by Schwarze (cf.~[Kr-T, 9.1]). We present the formulae for 
the image of $b$ under $\Phi$ in $I^2(K)$ (as in [\TccII, (9.5)]).

$$(\Phi_*b)(g,h)=
\cases{
[w']-[t]-[t']+[tt'w']&
if $g=g_2(u,t,v)$, $h=g_2(u',t',v')$, $w'=-(v+u')\ne0$,\cr  
-[1]-[t]-[t']-[tt']&
if $g=g_2(u,t,v)$, $h=g_2(u',t',v')$, $w'=-(v+u')=0$,\cr
\pf(t,t')&
if $g=g_1(u,t)$, $h=g_2(u',t',v')$,\cr
\pf(t,t')&
if $g=g_2(u,t,v)$, $h=g_1(u',t')$,\cr
\pf(t,t')&
if $g=g_1(u,t)$, $h=g_1(u',t')$.
}\leqno(\a.\f)
$$
\medskip

Step 3:

We work in the extension $E=\MG/\com(\G)\to\G$. 
Let $H$ be the image of $h(-1)\in\MG$ in $E$.
The subgroup $\pfund/\G^{[,]}$ of $E$  is generated by its index-2 subgroup
$$\fund/\G^{[,]}\simeq \Eq(\G)\simeq I^2(K)$$ and the element $H$; we have
$2H=\pf(-1,-1)$, a consequence of $h(-1)^2=\{-1,-1\}$.

The extension $E\to\G$ has a preferred lift: the lift of $\MG\to\G$
given by (\a.8), (\a.9) composed with the projection map $\MG\to E$.
The cocycle corresponding to this lift is $\Phi_*b$. Observe that the lift
of $-I$ in $E$ is $H$: indeed, $-I=g_1(0,-1)$ lifts to $h(-1)$ in $\MG$,
and $h(-1)$ maps to $H$ in $E$.

Now we would like to find commutators of lifts to $E$ of anti--commuting pairs in $\G$.
All such pairs are described by Corollary 3.2. Due to Proposition \a.2 we only need to consider
pairs given by (3.5), (3.6).
Assume first that $\gamma,\delta\ne0$ (necessarily true for Stufe 2).
We have
$$\eqalign{
b&=b'=\pmatrix{0&t\cr -t^{-1}&0}=w(t)=g_2(0,t,0),\cr
a&=a'=\pmatrix{\gamma&t^2\delta\cr \delta&-\gamma}=
g_2(\gamma/\delta,-\delta^{-1},-\gamma/\delta),\cr
ba&=\pmatrix{t\delta& -t\gamma\cr -t^{-1}\gamma& -t\delta}=
g_2(-t^2\delta/\gamma,t/\gamma,t^2\delta/\gamma),\cr
-I&=\pmatrix{-1&0\cr 0&-1}=g_1(0,-1).}
\leqno(\a.\f)$$
Using (\a.10) we compute
$$\eqalign{
\Phi_*b(a,b)&=[\gamma/\delta]-[-\delta^{-1}]-[t]+[-\gamma\delta^{-2}t]=
[\gamma\delta]+[\delta]-[t]-[\gamma t];\cr
\Phi_*b(b,a)&=[-\gamma/\delta]-[t]-[-\delta^{-1}]+[\gamma\delta^{-2}t]=
-[\gamma\delta]-[t]+[\delta]+[\gamma t];\cr
\Phi_*b(-I,ba)&=\pf(-1,t/\gamma)=[1]+[1]-[\gamma t]-[\gamma t].}
\leqno(\a.\f)
$$
Plugging this into (\a.1) we get
$$[\ov{a},\ov{b}]=H-2[1]+2[\gamma t]+[\gamma\delta]+[\delta]-[t]-[\gamma t]
+[\gamma\delta]+[t]-[\delta]-[\gamma t]=H-2[1]+2[\gamma\delta].
\leqno(\a.\f)$$
So, to pass from $\pfund/\com(\G)$ to $\Eq(\PG)$ we have to quotient by
the relations
$$H-2[1]+2[\gamma\delta]=0,\leqno(\a.\f)$$
where $(\gamma, \delta)$ run through all possible non-zero solutions to
$$\gamma^2+(t\delta)^2=-1\leqno(\a.\f)$$
with arbitrary non-zero $t$ (which can also vary).
The equation $x^2+y^2=-1$ is solvable in $K$ because $\stufe(K)\le2$;
consequently, it has infinitely many solutions, and so does (\a.15).
Now, once we have one non-zero solution of (\a.15),
we can change to another one
with an arbitrary non-zero value of $\delta$, just by
$(\gamma,\delta,t) \mapsto(\gamma,\delta s,t/s)$.
Therefore $\gamma\delta$ in (\a.14)
can be an arbitrary non-zero scalar:
$$H-2[1]+2[\eta]=0,\qquad \eta\in\dotK.\leqno(\a.\f)$$
Putting $\eta=1$ we get $H=0$. Taking this into account, and noting that
$2[1]-2[\eta]=\pf(\eta,\eta)$ we get the final result in the case $\stufe(K)=2$.

We state the result of the above discussion now, address an obvious question that it raises in a remark,
and return to the remaining case $\stufe(K)=1$
afterwards.

\smallskip
\bf Theorem \a.\t\sl\par
Let $K$ be an infinite field of characteristic $\ne2$ and Stufe 2. Then
$$\Eq(\PG)\simeq I^2_+(K)/\langle H,\pf(\eta,\eta)\mid \eta\in\dotK\rangle
\simeq I^2(K)/\langle\, \pf(\eta,\eta)\mid \eta\in\dotK\rangle
=I^2(K)/2I(K).\leqno(\a.\f)$$
\rm\par
\smallskip
\bf Remark \a.\t. \rm \par
Fields $K$ that satisfy $I^2(K)=2I(K)$ are called \it 1-stable \rm
([EL2, Definition 3.8]). Not all fields of Stufe 2 are 1-stable as the following example shows.

Let $K=\Q(\epsilon)$, where $\epsilon={1\over2}(-1+i\sqrt{3})$.
Note that $K$ is a field of Stufe 2: $\pm i\not\in K$ and $\epsilon^2+(\epsilon^2)^2=-1$.
We will show that $q=\pf(2+\epsilon,5)$ is non-zero in  $I^2(K)/2I(K)$.

We begin by checking that both $q=\langle 1,-(2+\epsilon),-5,10+5\epsilon\rangle$
and $q'=\langle -1,-(2+\epsilon),-5,10+5\epsilon\rangle$ are anisotropic over $K$; 
to do this we show that they are anisotropic over the 5-adic completion $K_5$ 
(note that 5 is prime in $K$). The residue field
${\cal O}_K/5{\cal O}_K$
is a quadratic extension ${\bf F}_{25}={\bf F}_5(\epsilon)$ of ${\bf F}_5$.
By [Lam, Proposition VI.1.9] we need to check that
$\langle \pm 1,-(2+\epsilon)\rangle$ and 
$\langle -1,2+\epsilon\rangle$ are anisotropic over ${\bf F}_{25}$.
These conditions are equivalent to the fact that $\pm(2+\epsilon)\not\in\dot{\bf F}_{25}^2$.
One checks directly that $-(2+\epsilon)$ is a generator of the cyclic group $\dot{\bf F}_{25}$,
hence is not a square. Next, $-1$ is a square in ${\bf F}_5$ (hence also in ${\bf F}_{25}$),
therefore $2+\epsilon=-(2+\epsilon)\cdot (-1)$ is not a square in ${\bf F}_{25}$.
This completes the proof that $q$ and $q'$ are anisotropic over $K$.

Now the argument is as follows. If $q$ were an element of $2I(K)\subseteq W(K)\cdot\langle 1,1\rangle$,
then, by [EL1, Theorem 1.4], $q$ could be decomposed as $\langle 1,1\rangle\cdot \sigma$ with $\sigma=\langle 1,b\rangle$
for some $b$. Then $q=\langle 1,1, b,b\rangle$; by Witt cancellation we deduce
$\langle -(2+\epsilon),-5,10+5\epsilon\rangle=\langle 1,b,b\rangle$, hence
$q'=\langle -1,-(2+\epsilon),-5,10+5\epsilon\rangle=
\langle -1,1,b,b\rangle$.
The latter form is, however, isotropic, while $q'$ is not. Contradiction.

\medskip
This concludes our discussion of $\Eq(\PG)$ when $\stufe(K)=2$.
Determining this group in the remaining case $\stufe(K)=1$ requires an extra calculation. In that case
step 3 of our procedure is not yet finished.
We have to take into account all the relations found above in the case
$\stufe(K)=2$---but if $\stufe(K)=1$, i.e.~$-1\in\dotK^2$, they only give $H=0$, because
$$\pf(\eta,\eta)=[1]-[\eta]-[\eta]+[1]=[1]-[\eta]+[-\eta]-[-1]
=[1]-[\eta]+[\eta]-[1]=0.\leqno(\a.\f)$$  
There are, however, additional relations, coming from the solutions
of (\a.15) with $\gamma=0$ or with $\delta=0$. In those cases formulae (3.5), (3.6) give:

$$\eqalign{
b&=b'=\pmatrix{0&t\cr -t^{-1}&0}=w(t)=g_2(0,t,0),\cr
a&=a'=\pmatrix{\gamma&0\cr 0&-\gamma}=g_1(0,\gamma)
\quad\hbox{\rm or}\quad
a=a'=\pmatrix{0&t^2\delta\cr \delta&0}=
g_2(0,-\delta^{-1},0),\cr
ba&=\pmatrix{0& -t\gamma\cr -t^{-1}\gamma& 0}=
g_2(0,t/\gamma,0),
\quad\hbox{\rm or}\quad
\pmatrix{t\delta& 0\cr 0& -t\delta}=
g_1(0,t\delta),
\cr
-I&=\pmatrix{-1&0\cr 0&-1}=g_1(0,-1).}
\leqno(\a.\f)$$

Using (\a.10) we compute

$$\eqalign{
\Phi_*b(a,b)&=\pf(\gamma,t)
\quad\hbox{\rm or}\quad
\pf(\delta^{-1},-t)=\pf(\delta,t);\cr
\Phi_*b(b,a)&=\pf(t,\gamma)
\quad\hbox{\rm or}\quad
\pf(-t,\delta^{-1})=\pf(t,\delta);\cr
\Phi_*b(-I,ba)&=\pf(-1,t/\gamma)=0
\quad\hbox{\rm or}\quad
\pf(-1,t\delta)=0.}
\leqno(\a.\f)
$$
Plugging this into (\a.1) we get
$$
[\ov{a},\ov{b}]=H-\Phi_*b(-I,ba)+\Phi_*b(a,b)-\Phi_*b(b,a)=H.\leqno(\a.\f)$$
Thus, the additional relations all reduce to $H=0$. We have proved:
\smallskip
\bf Theorem \a.\t.\sl\par
Let $K$ be an infinite field of characteristic $\ne2$ and Stufe 1. Then
$$\Eq(\PG)\simeq I^2(K).\leqno(\a.\f)$$
The isomorphism is induced by the epimorphism
$\pfund\to I^2(K)$ that sends $\{s,t\}$ to $\pf(s,t)$ and $h(-1)$ to $0$.
\rm\par
Note that the epimorphism in Theorem \a.5 is compatible with the one from
Corollary 3.6. Indeed, that map sends $h(-1)$ to $\langle 1,1\rangle$;
for Stufe 1 we have $\langle 1,1\rangle=0$ and $I^2_+(K)=I^2(K)$.

\bigskip

\medskip
\bf 5. The Witt class for $PSL(2,K)$.\rm
\medskip
\def\a{5}
\m=1
\n=1

In this section we investigate the relationship between the universal
equicommutative class of $\PG$ and the Witt class. Recall that
for $\G$ these classes are equivalent (cf [\TccII, Theorem A]).

The Witt cocycle on $\PG$ can be defined by
$$w(\pm A, \pm B)=[-A_{21}(AB)_{21}B_{21}].\leqno(\a.\f)$$
(By $M_{21}$ we denote the 21-entry of the matrix $M$.) 
The same formula defines the (pulled-back) cocycle $w$ on $\G$.
\smallskip
\bf Proposition \a.\t.\sl\par
The Witt cocycle on $\PG$ is equicommutative if and only if $\stufe(K)\ne2$.
\rm\par
Proof. If $\stufe(K)>2$ then (see Lemma 3.1) each commuting pair $(\pm A, \pm B)$
in $\PG$ comes from a commuting pair $(A,B)$ in $\G$. Using (\a.1) we can see
that $w(A,B)=w(B,A)$.

If $\stufe(K)=1$ then $[t]=[-t]$ for all $t\in\dotK$. Now anti-commuting pairs
$(A,B)$ in $\G$ exist, but for such pair we get from (\a.1):
$$w(A,B)=[-A_{21}(AB)_{21}B_{21}]=[A_{21}(BA)_{21}B_{21}]=w(B,A).\leqno(\a.\f)$$

If $\stufe(K)=2$ then there exist an anti-commuting pair $(A,B)$, we have
$$w(A,B)=[t],\quad w(B,A)=[-t]\leqno(\a.\f)$$
for some $t\in\dotK$ (see (3.5),(3.6); we have $\delta,\gamma\ne0$ there,
as otherwise $\stufe(K)=1$). We have $[t]\ne[-t]$---otherwise $-1\in\dotK^2$.
\qed
\medskip

One may guess that the Witt class and the universal equicommutative
cohomology class
of $\PG$ are equivalent if $\stufe(K)\ne2$. One may also be mildly puzzled by the relation
between these two classes for fields of Stufe $2$. 
Now we proceed to confirm the guess and solve the puzzle.
\medskip
Even for $\G$, the equivalence of the Witt class (represented by the Witt cocycle (\a.1))
and the universal equicommutative class (represented by $\Phi_*b$, see (4.10))
is not so straightforward. The coefficient groups of theses two cohomology classes
are different ($W(K)$ and $I(K)$). Compatibility is achieved by adding to the Witt
cocycle the so-called
Nekov\'a$\check{\rm r}$ correction: the coboundary $\delta n$,
where $n$ is given by
$$n(g)=\cases{
[|e,ge|] & if $|e,ge|\ne0$,\cr
[1]-[t]& if $ge=te$ for some $t\in\dotK$. 
}\leqno(\a.\f)$$
Here $e={1\choose0}$.
One easily checks that
$$n(g_1(u,t))=[1]-[t],\qquad n(g_2(u,t,v))=[-t^{-1}]=-[t].\leqno(\a.\f)$$
The following proposition (and its proof) can be found in [Kr-T, Theorem 23] or in
[\TccII, Proposition 9.2]; the idea goes back to [Ne] and [Ba].
\smallskip
\bf Proposition \a.\t.\sl\par
On $\G$ we have the cocycle equality $w+\delta n=\Phi_*b$. 
\rm\par
\smallskip
One consequence of this proposition is that the Witt class $[w]\in H^2(\G,W(K))$ has a reduction
to $H^2(\G,I^2(K))$: the class $[w+\delta n]=[\Phi_*b]\in H^2(\G,I^2(K))$ maps to $[w]$ under the inclusion
of cohomology groups induced by the inclusion $I^2(K)\to W(K)$ of coefficient groups
(cf [\TccII, Proposition 9.3, Definition 9.4]). 
\smallskip
We would like to prove a similar result for $\PG$. To this end, we adapt the cochain $n$ to $\PG$.
Define $\wt{n}\colon\PG\to W(K)$ by $\wt{n}(\pm A)=n(\sigma(A)A)$.
We have
$$
\eqalign{
(w+\delta\wt{n})(\pm A, \pm B)&=
w(\pm A, \pm B)+\wt{n}(\pm A)-\wt{n}(\pm AB)+\wt{n}(\pm B)\cr
&=w(\sigma(A)A, \sigma(A)B)+{n}(\sigma(A) A)-{n}(\sigma(AB) AB)+{n}(\sigma(B) B)\cr
&=(w+\delta n)(\sigma(A)A, \sigma(A)B)+{n}(\sigma(A)\sigma(B) AB)-{n}(\sigma(AB) AB)\cr
&=\Phi_*b(\sigma(A)A, \sigma(A)B)+{n}(\sigma(A)\sigma(B) AB)-{n}(\sigma(AB) AB).
}\leqno(\a.\f)
$$
We define the Moore cocycle $b_{\PG}$ on $\PG$ as the cocycle of the extension
$\MG\to\PG$ corresponding to the lift
$\pm A\mapsto \ov{\sigma(A)A}$. Let us express $b_{\PG}$ in terms of the usual Moore cocycle
$b$ on $\G$:
$$
\eqalign{
b&_{\PG}(\pm A,\pm B)=\ov{\sigma(A)A}\cdot \ov{\sigma(B)B}\cdot\ov{\sigma(AB)AB}^{-1}\cr
&=\ov{\sigma(A)A}\cdot \ov{\sigma(B)B}\cdot\ov{\sigma(A)A\sigma(B)B}^{-1}
\cdot\ov{\sigma(A)A\sigma(B)B}\cdot\ov{\sigma(AB)AB}^{-1}\cr
&=b(\sigma(A)A,\sigma(B)B)
\cdot\ov{\sigma(A)A\sigma(B)B(\sigma(AB)AB)^{-1}}\cdot
b(\sigma(A)A\sigma(B)B(\sigma(AB)AB)^{-1},\sigma(AB)AB)^{-1}\cr
&=b(\sigma(A)A,\sigma(B)B)
\cdot\ov{\sigma(A)\sigma(B)\sigma(AB)I}\cdot
b(\sigma(A)\sigma(B)\sigma(AB)I,\sigma(AB)AB)^{-1}.
}\leqno(\a.\f)
$$
(We have used the general cocycle
formula $\ov{g}\cdot\ov{h}^{-1}=\ov{gh^{-1}}\cdot b(gh^{-1},h)^{-1}$.)

We know that $\ov{-I}=h(-1)$. We extend the map $\Phi\colon \fund\to I^2(K)$
to $\Psi\colon\pfund\to I^2_+(K)$ by $\Psi(h(-1))=\langle 1,1\rangle$.
This $\Psi$ induces a map $\pfund/\com(\G)\to I^2_+(K)$---an isomorphism for
Stufe $>1$, and an epimorphism with kernel $\Z/2=\{1,h(-1)\}$ for Stufe $1$.
\smallskip
\bf Theorem \a.\t\sl\par
On $\PG$ we have the cocycle equality
$$\Psi_*b_{\PG}=w+\delta\wt{n}.\leqno(\a.\f)$$
In particular, $w+\delta\wt{n}$ is an $I^2_+(K)$-valued cocycle,
providing the reduction of the Witt class of $\PG$ to $H^2(\PG,I^2_+(K))$.
\rm\par
Proof. For a pair $(\pm A, \pm B)$ in $\PG$  that satisfies
$\sigma(AB)=\sigma(A)\sigma(B)$ this follows
from (\a.6), (\a.7) (using $\ov{I}=1$ and $b(I,*)=1$).

Suppose that $\sigma(AB)=-\sigma(A)\sigma(B)$. Then, again due to (\a.6), (\a.7), we have to
check that
$$\Psi_*(\ov{-I})-\Psi_*b(-I,-\sigma(A)\sigma(B)AB)=n(\sigma(A)A\sigma(B)B)-n(-\sigma(A)A\sigma(B)B).
\leqno(\a.\f)$$
We consider two cases and compute the left- and right-hand-side
of this formula using (4.10) and (\a.5).
\smallskip
Case 1. $\sigma(A)A\sigma(B)B=g_1(u,t)$. Then $-\sigma(A)A\sigma(B)B=-g_1(u,t)=g_1(u,-t)$.
Then
$$
\eqalign{
LHS&= \langle 1,1\rangle - \Phi_*b(g_1(0,-1),g_1(u,-t))=
2[1]-\pf(-1,-t)=2[1]-([1]-[-1]-[-t]+[t])=-2[t],\cr
RHS&=n(g_1(u,t))-n(g_1(u,-t))=[1]-[t]-([1]-[-t])=-2[t].
}
$$

Case 2. $\sigma(A)A\sigma(B)B=g_2(u,t,v)$. Then $-\sigma(A)A\sigma(B)B=-g_2(u,t,v)=g_2(u,-t,v)$.
Then
$$
\eqalign{
LHS&= \langle 1,1\rangle - \Phi_*b(g_1(0,-1),g_2(u,-t,v))=
2[1]-\pf(-1,-t)=2[1]-([1]-[-1]-[-t]+[t])=-2[t],\cr
RHS&=n(g_2(u,t,v))-n(g_2(u,-t,v))=-[t]-(-[-t])=-2[t].
}
$$

\qed
\medskip
The next theorem is important: it explains the relation between the Witt class and the universal
equicommutative class.
\smallskip
\bf Theorem \a.\t\sl\par
The Witt class, reduced to $I^2_+(K)$ as in Theorem \a.3:
\item{a)}{is equal to the image of the universal Moore class of $\PG$
under the coefficient map $$\Psi\colon\pfund\mathop{\longrightarrow} \pfund/\com(\G)\mathop{\longrightarrow}
I^2_+(K);$$}
\item{b)}{if $\stufe(K)\ne2$---is equivalent to the universal equicommutative class of $\PG$;}
\item{c)}{if $\stufe(K)=2$---is not equicommutative, but its image in
$$
I^2_+(K)/\langle H,\pf(\eta,\eta)\mid \eta\in\dotK\rangle
\simeq
I^2(K)/2I(K)$$
under the coefficient quotient map is the universal equicommutative class
of $\PG$.}
\rm\par
Proof: Part a) follows from Theorem \a.3. Then part b) follows from Corollary 3.6 and
Theorem 4.5, and part c) follows from Theorem 4.3. 
\qed
\smallskip

\bf Remark \a.\t.\rm\par
We shall make explicit a universal property of the Witt class for
$K$ of Stufe 2. Call a pair $(\pm A,\pm B)$ of elements
of $\PG$ \it $+$-commuting, \rm if $AB=BA$ (this property does not
depend on the choice of the lifts to $\G$). Formula (3.2)
expresses $\com(\G)$ in terms of $+$-commuting pairs in $\PG$.
One may call a cocycle $c\in H^2(\PG,A)$ \it $+$-equicommutative, \rm
if for each $+$-commuting pair $(\pm A, \pm B)$ one has
$c(\pm A, \pm B)=c(\pm B, \pm A)$; a cohomology class would be
\it $+$-equicommutative, \rm if it is represented by a $+$-equicommutative
cocycle. Then Theorem 5.4.a) implies that for fields of Stufe 2
the Witt class is universal $+$-equicommutative: it is $+$-equicommutative,
and for every $+$-equicommutative class $\gamma\in H^2(\PG,A)$ 
there exists a unique coefficient homomorphism $I^2_+(K)\to A$ that maps
the Witt class to $\gamma$. Formally, this universality holds also
when $\stufe(K)>2$ (as the properties of $+$-commuting and commuting
are then the same), but not for $\stufe(K)=1$ (where these properties differ,
and the Witt class is universal equicommutative, $+$-equicommutative,
but not universal $+$-equicommutative).

\bigskip

\medskip
\bf 6. Milnor--Wood bounds on the Witt class. \rm
\medskip
\def\a{6}
\m=1
\n=1

Suppose that we are given a class $u\in H^2(G,U)$.
Then, for any representation $\rho\colon\pi_1\Sigma\to G$
of the fundamental group of an oriented closed surface $\Sigma$, we can
define
$$\rho^*u\in H^2(\pi_1\Sigma,U)\simeq H^2(\Sigma,U)\simeq U.\leqno(\a.\f)$$
The first isomorphism reflects the fact that $\Sigma=K(\pi_1\Sigma,1)$;
the second is the evaluation on the fundamental class $[\Sigma]\in H_2(\Sigma,\Z)$.
We will use notation $u(\rho)$ for the element $\rho^*u$ pushed to $U$
by these isomorphisms. A natural question is: what is the set of possible values
of $u(\rho)$ for a fixed $\Sigma$ and varying $\rho$?

This question can be rephrased in terms of (flat, principal) $G$-bundles.
To any representation $\rho\colon\pi_1\Sigma\to G$
one can associate a $G$-bundle $\xi_\rho$ over $\Sigma$.
One obtains a bijective correspondence between classes of such representations
modulo conjugation in $G$ and isomorphism classes of $G$-bundles over $\Sigma$.
The invariant $u(\rho)\in U$ does not change under $G$-conjugation
(internal automorphisms act trivially on cohomology with constant coefficients),
therefore one can set $u(\xi_\rho):=u(\rho)$, defining a $U$-valued invariant
of $G$-bundles over $\Sigma$.

The question can now be stated as follows: what is the range of $u(\xi)$ for
$G$-bundles over a closed surface of genus $g$? The most classical instance
of this problem is the Milnor--Wood inequality (cf [Mil]).
In [\TccII, Theorem 11.6] a partial answer is provided for the Witt class
of $\G$. We will prove an analogous theorem for $\PG$. To state it we need
a norm on our $U=I^2_+(K)$ that will measure how large the invariant $u(\xi)$
can get when the genus $g$ grows.

The Witt group $W(K)$ has a natural norm. That norm has two equivalent descriptions:
as the word norm with respect to the standard system of generators $([a]\mid a\in\dotK)$,
or as the dimension of an anisotropic representative (cf [\TccII, Lemma 11.5]).
We shall equip the subgroups $I^2(K)$ and $I^2_+(K)$ with the restrictions of this norm.
Since $I^2(K), I^2_+(K)\subseteq I(K)$, these restrictions take only even values.
We can now state the promised:
\smallskip
\bf Theorem \a.\t. \sl\par
Let $K$ be an infinite field.
\item{(a)} The Witt class of any flat $PSL(2,K)$-bundle over an oriented
closed surface of genus $g$ has norm $\le 4(g-1)+2$.
\item{(b)} 
The set of Witt classes of flat $PSL(2,K)$-bundles over an oriented
closed surface of genus $g$ contains the set of elements
of $I^2_+(K)$ of norm $\le 4(g-1)$.
\smallskip\rm

The part of this theorem not covered by [\TccII, Theorem 11.6] is
the construction of bundles whose Witt classes realize elements of
$I^2_+(K)\setminus I^2(K)$. These are $\PG$-bundles that cannot admit
a reduction to $\G$. They will be constructed by splitting $\Sigma=S_1\cup S_2$
and gluing some $\G$-bundles $\xi_1$, $\xi_2$ over $S_1$, $S_2$, whose monodromies
along the common boundary $\partial S_1=\partial S_2$ will match in $\PG$ but not in $\G$.
The next lemma will be used to find the Witt class of bundles glued in this way;
it can---and will---be stated in a more general setting. The statement
requires some introductory definitions.

Suppose that $E\to G$ is a central extension with kernel $U$, $g\mapsto\ov{g}$
is a lift, and $c\in Z^2(G,U)$ is the associated cocycle (as in (1.4)).
We assume that $\ov{1}=1$; then the cocycle is normalized:
$c(1,g)=c(g,1)=1$. Now let $S$ be a compact oriented surface with
one (compatibly oriented) boundary component $\partial S$ and a chosen
basepoint $s\in\partial S$. Let $\xi$ be a $G$-bundle over $S$; let
a fibre trivialisation $\xi_s\simeq G$ be chosen. We shall associate
to these data an element $\ov{c}^G(\xi)\in U$, called a
\it relative class with standard framing \rm (cf [TccII, Definition 4.1]).
We choose in $S$ a \it standard loop collection: \rm a collection
of loops $(x_i,y_i\mid i=1,\ldots,g)$, based at $s$, cutting
$S$ into a $(4g+1)$-gon, generating $\pi_1(S)$, satisfying $\prod_i[x_i,y_i]=\partial S$
in $\pi_1(S)$. We denote by $X_i,Y_i\in G$ the monodromies along these loops,
and by $W$ the monodromy along $\partial S$. We then put:
$$\ov{c}^G(\xi)=\left(\prod_i[\ov{X_i},\ov{Y_i}]\right){\ov{W}}^{-1}\in U.\leqno(\a.\f)$$

Now suppose that $G$ itself is a central extension of another group $PG$:
an epimorphism $G\ni g\mapsto [g]\in PG$ is given, with kernel 
central in $G$
(our prime example is $\G\to\PG$ with kernel $\{I,-I\}$). Suppose also that
the cocycle $c$ is the pull-back to $G$ of a normalized cocycle $c^{PG}\in Z^2(PG,U)$,
so that $c(g,h)=c^{PG}([g],[h])$ (just as it happens for the Witt cocycle (5.1) in our prime example).
Let $PE\to PG$ be a central extension with kernel $U$ and a lift $[g]\mapsto\ov{[g]}$,
corresponding to the cocycle $c^{PG}$, ie such that
$$\ov{[g]}\cdot\ov{[h]}=\ov{[gh]}\cdot c^{PG}([g],[h]).\leqno(\a.\f)$$
(A standard way to define such an extension is to put $PE=PG\times U$
with multiplication
$$([g],u)\cdot([g'],u'):=([g][g'], uu'c^{PG}([g],[g']))\leqno(\a.\f)$$
and with the first--factor projection $PG\times U\to PG$ as the extension epimorphism.)
Every element of $E$ is uniquely expressible as $\ov{g}u$ for some
$g\in G$, $u\in U$. We define $\pi\colon E\to PE$ by
$$\pi(\ov{g}u)=\ov{[g]}u.\leqno(\a.\f)$$
Straightforward calculation shows that this map
is an epimorphism. 
Clearly, $\pi|_U={\rm Id}_U$.

Denote 
by $\tau_c^{PG}\in H^2(PG,U)$ the cohomology class of $c^{PG}$.

\smallskip
\bf Lemma \a.\t.\sl\par
Let $\xi_1$, $\xi_2$ be $G$-bundles over oriented surfaces $S_1$, $S_2$ with boundary
monodromies that have the same image in $PG$; let $\xi_1^{PG}$, $\xi_2^{PG}$ be the induced $PG$-bundles.
Let $\Sigma=S_1\cup_{\partial} S_2$ be oriented compatibly
with $S_1$ and opposite to $S_2$. Let $\xi=\xi_1^{PG}\cup\xi_2^{PG}$ be the glued $PG$-bundle
over $\Sigma$. Then
$$\tau_c^{PG}(\xi)=\ov{c}^G(\xi_1)\,\ov{c}^G(\xi_2)^{-1}.\leqno(\a.\f)$$
\rm
\smallskip
Proof. Let $(x_i,y_i)$ be a standard loop collection for $S_1$, and let $(X_i,Y_i)$ be
the corresponding collection of monodromies of $\xi_1$. Let $W$ be the boundary monodromy
of $\xi_1$.
Using $\ov{c}^G(\xi_1)\in U$, $\pi|_U={\rm Id}_U$, and the fact that $\pi$ is a homomorphism, we calculate:
$$\eqalign{
\ov{c}^G(\xi_1)=\pi(\ov{c}^G(\xi_1))&=
\pi\left(\left(\prod_i\left[\ov{X_i},\ov{Y_i}\right]\right)\cdot \ov{W}^{-1}\right)\cr
&=\left(\prod_i\left[\pi\left(\ov{X_i}\right),\pi\left(\ov{Y_i}\right)\right]\right)\cdot \pi\left(\ov{W}\right)^{-1}
=\left(\prod_i\left[\ov{[X_i]},\ov{[Y_i]}\right]\right)\cdot \ov{[W]}^{-1}=\ov{c}^{PG}(\xi_1^{PG}).}
\leqno(\a.\f)$$
Similarly, $\ov{c}^G(\xi_2)=\ov{c}^{PG}(\xi_2)$.
Now [\TccII, Lemma 4.4] gives:
$$\tau_c^{PG}(\xi)=\ov{c}^{PG}(\xi_1)\ov{c}^{PG}(\xi_2)^{-1}=\ov{c}^G(\xi_1)\ov{c}^G(\xi_2)^{-1}.\leqno(\a.\f)$$
\qed

\medskip

Proof of Theorem \a.1.

We denote by $w^{\G}$, $w^{\PG}$ the Witt cocycles, and by
$\w^{\G}$, $\w^{\PG}$ the corresponding cohomology classes.
Since $w^{\G}$ is the pull-back of $w^{\PG}$, we have
$\w^{\G}(\xi)=\w^{\PG}(\xi^{\PG})$ for all $\G$-bundles $\xi$.
\smallskip

Part (a) is proved by an easy classical argument, exactly as in
[\TccII, Theorem 11.6].
\smallskip
To prove part (b), consider an element $q\in I^2_+(K)$ of norm $\le 4(g-1)$.
We claim that $q=\w^{\PG}(\xi)$ for some $\PG$-bundle $\xi$ over a surface of genus $g$. 

If $q\in I^2(K)$, then by [\TccII, Theorem 11.6] we know that
$q=\w^{\G}(\xi)$ for some $\G$-bundle $\xi$ over a surface of genus $g$;
but $\w^{\G}(\xi)=\w^{\PG}(\xi^{\PG})$, hence the claim follows in this case.

Suppose then that $q\in I^2_+(K)\setminus I^2(K)$.
We follow the proof of [\TccII, Theorem 11.6].
Let
$$q=\sum_{i=1}^g([\alpha_i]+[\beta_i])+\sum_{j=2}^{g-1}([\gamma_j]+[\delta_j]).\leqno(\a.\f)$$
Since $q\in I^2_+(K)\setminus I^2(K)$, we have $d_\pm(q)=-1$
(indeed, $q=\langle 1,1\rangle+q'$ for some $q'\in I^2(K)$; hence
$d_\pm(q)=d_\pm(\langle 1,1\rangle)d_\pm(q')=(-1)\cdot 1=-1$). Therefore
$$-1=d_\pm(q)=\prod_i(\alpha_i\beta_i)\prod_j(\gamma_j\delta_j).\leqno(\a.\f)$$
We construct $\G$-bundles $\zeta_{g-1}$, $\xi_g$ exactly as in the proof of
[\TccII, Theorem 11.6].
Their relative Witt classes are
$$
\ov{w}^{\G}(\zeta_{g-1})=\sum_{i=1}^{g-1}([\alpha_i]+[\beta_i])+\sum_{j=2}^{g-1}([\gamma_j]+[\delta_j]),\qquad
\ov{w}^{\G}(\xi_g)=[-\alpha_g]+[-\beta_g].
\leqno(\a.\f)$$
Their boundary monodromies 
are diagonal matrices 
with entries $u_{g-1}$, $z_g$, respectively, whose classes in $\dotK/\dotKK$
satisfy:
$$[u_{g-1}]=[-\prod_{i=1}^{g-1}(\alpha_i\beta_i)\prod_{j=1}^{g-1}(\gamma_j\delta_j)],\qquad
[z_g]=[-\alpha_g\beta_g].\leqno(\a.\f)$$
Moreover (cf.~[\TccII, Lemma 11.2, Corollary 11.4]), the exact values of $u_{g-1}$, $z_g$
can be chosen almost arbitrarily within the appropriate class of $\dotK/\dotKK$
(there are finitely many exceptional values). Observe that
$$[u_{g-1}z_g]=[\prod_{i=1}^g(\alpha_i\beta_i)\prod_{j=2}^{g-1}(\gamma_j\delta_j)]=[-1].\leqno(\a.\f)$$
Therefore, we may choose $u_{g-1}$, $z_g$ so that $z_g=-u_{g-1}$; then $\zeta_{g-1}^{\PG}$ and $\xi_g^{\PG}$
glue to a $\PG$-bundle $\xi$. We use Lemma \a.2:
$$
\eqalign{
\w{}^{\PG}(\xi)&=\ov{w}^{\G}(\zeta_{g-1})-\ov{w}^{\G}(\xi_g)\cr
&=\sum_{i=1}^{g-1}([\alpha_i]+[\beta_i])+\sum_{j=2}^{g-1}([\gamma_j]+[\delta_j])
-([-\alpha_g]+[-\beta_g])=q}
\leqno(\a.\f)$$
\qed
\smallskip
The discrepancy between the bounds in part a) (``the weak Milnor--Wood inequality'') and part b) (``the strong Milnor--Wood inequality'') of Theorem \a.1
leads to an obvious question: can elements of norm $4(g-1)+2$ be realized
as Witt classes of bundles over $\Sigma_g$? If $K=\R$
they cannot, as the classical Milnor--Wood inequality attests. The argument
is by stabilising the weak Milnor--Wood inequality, ie using it on a covering
space---we shall use this argument to clarify the situation for $K=\Q$.
Then, in Remark \a.5, we shall observe that the calculations in
Section 3 legally yield
a counterexample in genus $1$
(which raises more questions than it answers).
\smallskip
\bf Proposition \a.\t.\sl\par
The set of Witt classes of flat $PSL(2,\Q)$-bundles over an oriented
closed surface of genus $g$ is equal to the set of elements
of $I^2_+(\Q)$ of norm $\le 4(g-1)$.
\rm\par
\smallskip
Proof. Since $\stufe(\Q)=+\infty$, the Witt class is equicommutative
and all bundles over $T^2$ have Witt class $0$.

For $g>1$ we use Theorem \a.1a). We have to prove that there are
no bundles over $\Sigma_g$ with Witt class
of norm $4(g-1)+2$. But, for elements of $W(\Q)$ of norm $>4$ the norm
is equal to $\pm$ the signature $\sigma$ (as a consequence of Meyer's Theorem
[MH, Corollary II.3.2]). Now if we had a $q\in I^2_+(\Q)$ with norm
$\|q\|=4(g-1)+2\ge6$ realized as the Witt class of
a bundle $\xi$ over $\Sigma_g$,
then $2q$ would be equal to the Witt class of
$\widetilde{\xi}$---the pull-back of $\xi$ to a 2-fold covering
$\widetilde{\Sigma}$ of $\Sigma$. Then, on the one hand,
$\|2q\|=|\sigma(2q)|=2|\sigma(q)|=4(2g-2)+4$, while on the other hand,
the genus of $\widetilde{\Sigma}$ would be $2g-1$---altogether contradicting
the weak Milnor--Wood bound from Theorem \a.1.a).\qed
\smallskip

\smallskip
\bf Remark \a.\t.\rm\par
The last part of the proof of Proposition \a.3 shows
the following general statement: Let $K$ be an ordered field;
then the signature of the Witt class satisfies the strong
Milnor--Wood inequality.

\smallskip
\bf Remark \a.\t.\rm\par
Suppose that $\stufe(K)=2$ and that the genus of $\Sigma$
is $1$ ($\Sigma=T^2$). In this case, we shall prove that each 
element of norm $2$ in $I^2_+(K)$ is the Witt class of a $\PG$-bundle
over $T^2$.

Let $q\in I^2_+(K)$ have norm 2; then $q=\langle\eta,\mu\rangle$
for some $\eta,\mu\in\dotK$. By Hauptsatz ([Lam, Theorem X.5.1])
elements of $I^2(K)$ have norm at least 4; therefore
$q\in \langle 1,1\rangle+I^2(K)$,
hence $-\eta\mu=d_\pm(q)=-1$ (equality in $\dotK/\dotKK$). This implies that
$q=\langle\eta,\eta\rangle$ (equality in $W(K)$). However, in Section 3, Step 3
(esp.~(3.13), (3.16)) we have found, for each $\eta\in\dotK$,
a pair of anti-commuting matrices $a,b\in \G$ with
$[\ov{a},\ov{b}]=H-2[1]+2[\eta]=2[\eta]$ (equality in $I^2_+(K)$).
Such a pair $(a,b)$ is the monodromy pair of a $\PG$-bundle
over $T^2$ with Witt class $[\ov{a},\ov{b}]=2[\eta]$. 

\smallskip
\bf Remark \a.\t.\rm\par
The previous remark shows that the strong Milnor--Wood inequality
for the Witt class for $\PG$ fails for fields of Stufe $2$---a
counterexample is given over the surface
of genus $1$. We do not know whether there exist counterexamples
over surfaces of higher genus. 
It would also be of great interest to us to know whether the strong
Milnor--Wood inequality holds for the universal equicommutative
class (of $\PG$).

\bigskip
\bf References
\rm

\smallskip
\item{[Ba]}
J.~Barge,
Cocycle d'Euler et K2,
K-Theory {\bf 7} (1993), no.~1, 9--16.

\smallskip
\item{[Brown]}
 K.~Brown,
 Cohomology of groups.
 Graduate Texts in Mathematics, 87.
 Springer-Verlag, New York-Berlin, 1982.

\smallskip
\item{[\TccI]}
 J.~Dymara, T.~Januszkiewicz,
 Tautological characteristic classes I,
 Algebr.~Geom.~Topol. {\bf 24} (2024), no.~5, 2889--2932.

\smallskip
\item{[\TccII]}
 J.~Dymara, T.~Januszkiewicz,
 Tautological characteristic classes II: the Witt class,
 Math.~Ann. {\bf 390} (2024), 4463--4496.

\smallskip
\item{[EL1]}
 R.~Elman, T.Y.~Lam,
 Pfister Forms and K-Theory of Fields,
 J.~Algebra {\bf 23} (1972), 181--213.

\smallskip
\item{[EL2]}
 R.~Elman, T.Y.~Lam,
 Quadratic forms over formally real fields and pythagorean fields.
 Amer.~J.~Math. {\bf 94} (1972), 1155--1194.

\smallskip
\item{[Kr-T]}
 L.~Kramer, K.~Tent, A Maslov cocycle for unitary groups.
 Proc.~Lond.~Math.~Soc. (3) {\bf 100} (2010), no.~1, 91--115.

\smallskip
\item{[Lam]}
T.Y.~Lam,
Introduction to quadratic forms over fields.
Graduate Studies in Mathematics, 67. American Mathematical Society, Providence, RI, 2005.

\smallskip
\item{[Mil]}
 J.~Milnor,
 On the existence of a connection with curvature zero.
 Comment.~Math.~Helv. {\bf 32} 1958, 215--223.

\smallskip
\item{[MH]}
 J.~Milnor, D.~Husemoller,
 Symmetric bilinear forms. Ergebnisse der Mathematik und ihrer Grenzgebiete, Band 73. Springer--Verlag,
 New York--Heidelberg, 1973.

\smallskip
\item{[Moore]}
 C.~Moore, Group extensions of $p$-adic and adelic linear groups.
 Inst.~Hautes \'Etudes Sci.~Publ.~Math., {\bf 35} 1968, 157--222.

\smallskip
\item{[Ne]}
 J.~Nekov\'a$\check{\rm r}$,
 The Maslov index and Clifford algebras.
 (Russian) Funktsional.~Anal.~i Prilozhen. {\bf 24} (1990), no.~3, 36--44, 96;
 translation in Funct.~Anal.~Appl. {\bf 24} (1990), no.~3, 196--204 (1991). 

\smallskip
\item{[St]}
 R.~Steinberg, 
 Lectures on Chevalley groups. Notes prepared by John Faulkner and Robert Wilson.
 Revised and corrected edition of the 1968 original.
 University Lecture Series, 66. American Mathematical Society, Providence, RI, 2016. 

\smallskip
\item{[Su]}
 A.~Suslin,
 Torsion in $K_2$ of Fields,
 K-Theory {\bf 1} (1987), 5--29. 

\smallskip
\item{[Wei]}
 Ch.A.~Weibel,
 An introduction to homological algebra.
 Cambridge Stud.~Adv.~Math., 38.
 Cambridge University Press, Cambridge, 1994. xiv+450 pp.

\bye